\documentclass[12pt]{article}

\usepackage{amsmath}
\usepackage{amssymb}
\usepackage{amsthm}
\usepackage{bbm}
\usepackage[dvips]{graphicx}
\usepackage{color}
\usepackage{epsfig}
\usepackage[mathscr]{eucal}

\newtheorem{theorem}{Theorem}

\newtheorem{prop}[theorem]{Proposition}
\newtheorem{conj}{Conjecture}
\newtheorem{lemma}[theorem]{Lemma}
\newtheorem{cor}[theorem]{Corollary}

\newtheorem*{remark}{Remark}

\newcommand{\floor}[1]{\ensuremath{\left \lfloor {#1} \right \rfloor}}

\newcommand{\bv}{{\bf v}}

\newcommand{\bx}{{\bf x}}

\newcommand{\bL}{{\bf L}}

\newcommand{\cI}{{\mathcal I}}

\newcommand{\cT}{{\mathcal T}}
\newcommand{\cW}{{\mathcal W}}

\newcommand{\ignore}[1]{}
\newcommand{\pits}{\textsc{pits}}
\newcommand{\BE}{\mathbb{E}}
\newcommand{\BH}{\mathbb{H}}
\newcommand{\BP}{\mathbb{P}}

\allowdisplaybreaks

\title{Random Linear Extensions of Grids}
\author{Joshua N.\@ Cooper \\ \small ETH-Z\"urich} 

\begin{document}

\maketitle

\begin{abstract}
A grid poset -- or {\it grid} for short -- is a product of chains.  We ask, what does a random linear extension of a grid look like?  In particular, we show that the average ``jump number,'' i.e., the number of times that two consecutive elements in a linear extension are incomparable in the poset, is close to its maximum possible value.  The techniques employed rely on entropy arguments.  We finish with several interesting questions about this wide-open area.
\end{abstract}

\section{Introduction}

A {\it grid} poset -- or more succinctly, a {\it grid} -- is a product of chains.  Grids make appearances in a surprisingly diverse range of mathematics.  The lattice of submultisets of a multiset is a grid, including the extreme case of the Boolean lattice; the lattice of positive divisors of an integer is a grid; conjoint analysis (a branch of economic preference theory) is concerned with grids as sets of goods or products; Ferrers diagrams and plane partitions are lower-sets (i.e., order ideals) of two and three dimensional grids, respectively; the theory of poset order dimension (q.v. \cite{T92}) can be viewed as studying criteria for embeddability into grids.

In many of these contexts, one of the most important objects of study is the linear extension, i.e., a linear ordering of the vertices of the grid which respects the original ordering.  Linear extensions of posets are an object of extensive study in their own right, and the specialization to grids has raised many an interesting question.  Brightwell and Tetali \cite{BK01} recently pinned down the asymptotics for the number of linear extensions of the Boolean lattice, and the question of when such a linear extension is {\it representable} is a central issue in the area(s) known by the names ``qualitative'', ``subjective'', and ``comparative'' probability (q.v. \cite{F86,F96,F97,KPS59}).  Marketing researchers (e.g., \cite{KR76}) are interested in eliciting customers' personal extensions of product grids with few queries in order to determine pricing.  Linear extensions of a two dimensional grid correspond to Young Tableaux, objects to which a huge amount of attention has been devoted (and for which the interested reader is directed to \cite{S97} as a starting point).  The list goes on.

Frequently, it is a natural question to ask what a {\it random} linear extension of a grid looks like.  While specific aspects of this inquiry have been addressed in the literature (witness, e.g., \cite{B92, PR06}), the general setting, per se, does not appear to have been addressed, and so we hope to initiate some interest in the topic with the present work.  Naturally, there are many ways to attack the issue of describing a random structure, and we focus on one: the so-called {\it jump number}.

The {\it jump number} of a poset is the minimum number of times that a linear extension ``jumps'' from an element of the poset to an incomparable one.  (Formal definitions appear in the next section.)  This quantity, which in some sense measures how vertical or stratified a poset is, is well-studied.  (See, for example, \cite{CH79, N00, R96b}.)  We are interested in the question of what the jump number of a random linear extension of a grid is.  In particular, we show that, for a grid whose dimension is not too high, the jump number is almost surely close to the maximum possible number of jumps in any linear extension.  Notably, the proof uses information-theoretic arguments, an approach that appears to be gaining recognition as a powerful tool for attacking questions about linear extensions.

In the next section, we provide the necessary definitions.  The following section addresses the issue of counting linear extensions of grids, and the following section looks at the extreme jump numbers of such extensions.  Section 4 contains the main results on random linear extensions of grids, and we conclude with several open-ended questions in this area.

\section{Definitions}

For a positive integer $n$, the ``weak Bruhat order'' is the poset whose vertices are the permutations of $[n]=\{1,\ldots,n\}$, and so that $\pi_1 \succ \pi_2$ when there exists an $a \in [n-1]$ so that (1) $\pi_1(x) = \pi_2(x)$ whenever $x \neq a,a+1$, (2) $\pi_1(a) = \pi_2(a+1)$ and $\pi_1(a+1) = \pi_2(a)$, and (3) $\pi_2(a) < \pi_2(a+1)$.  The ``Bruhat graph'' $B_n$ is the Hasse diagram of the Bruhat order.

Suppose $P$ is a finite poset.  Define $\cT(P)$, the {\it transposition graph of $P$}, to be the induced subgraph of $B_n$ whose vertices are the linear extensions of $P$ for some (fixed) labelling of the vertices of $P$ by $[|P|]$.  Note that $\cT(P) = B_n$, the $1$-skeleton of the permutahedron, if $P$ is an $n$-element antichain.  We are interested in studying the linear extensions of $P = [m]^n$, a product of $n$ chains of length $m$, which we term a {\it grid}.

Let $r(\bx)$, for $\bx \in [m]^n$, denote the rank of $\bx = (x_1,\ldots,x_n)$, i.e., $1 + \sum_{j=1}^n (x_j - 1)$.  The $s^\textrm{th}$ rank-level $r^{-1}(s)$ we will denote by $\cW_s$.  We write $L^k$ for $L^{-1}(k)$, i.e., $k^\textrm{th}$ element of the ordering given by $L$, and $L^S = \bigcup_{j \in S} L^j$.

A {\it jump} of a linear extension $L$ of a poset is an integer $k \in [|P|]$, which we will refer to as a ``time'', so that $L^{k} \not \prec L^{k+1}$.  The {\it jump number} $s(P)$ of a poset $P$ is the least number of jumps over all linear extensions of $P$.  It is not hard to see that $s(P)$ is simply the minimum degree of $\cT(P)$.  The {\it average jump number} $\bar{s}(P)$, which corresponds to the average degree of the graph $\cT(P)$, is just the average number of jumps over all linear extensions of $P$.\\

\section{Vertex Count}

First, we give an estimate for the number of vertices in $\cT([m]^n)$.

\begin{prop} 
$$ \frac{1}{(ne)^{1/(n-1)}} \leq m^{-1} A(m,n)^{\frac{1}{(n-1)m^{n}}} \leq \frac{e}{2}.
$$
\end{prop}

\begin{proof} Note that a simple lower bound for the number of extensions of any ranked poset is the product of the factorials of its rank-level sizes.  Let $a_1,\ldots,a_{m(n-1)+1}$ denote the rank-level sizes of the poset $[m]^n$, i.e., $a_j = |\cW_j|$, the $j^\textrm{th}$ Whitney number.  We claim that, for any sequence of positive integers $b_1,\ldots,b_s$,
\begin{equation} \label{eq1}
\prod_{j=1}^s b_j! \geq \left ( \left ( s^{-1} \sum_{j=1}^s b_j \right ) ! \right)^s,
\end{equation}
where $x! = \Gamma(x+1)$ for $x \in \mathbb{R}$.  To see this, take the log of both sides:
$$
\sum_{j=1}^s \log(b_j!) \geq s \log \left ( s^{-1} \sum_{j=1}^s b_j \right ) !
$$
This equation clearly holds, by Jensen's inequality, if we can show that $\log x!$ is a convex function of $x$.  The second derivative of $\log \Gamma(x)$ is the {\it trigamma function}, with well-known representation
$$
\frac{d^2}{dx^2} \log \Gamma(x) = \sum_{j=0}^\infty \frac{1}{(j+x)^2},
$$
a quantity which is clearly positive for $x > -1$.  Therefore, we may conclude that (\ref{eq1}) holds.  By defining $a_{m(n-1)+2},\ldots,a_{mn} = 0$, we may apply this to the sequence $a_1,\ldots,a_{mn}$, yielding
$$
A(m,n) \geq \prod_{j=1}^{mn} a_j! \geq \left [ \left ( \frac{m^{n-1}}{n} \right ) ! \right]^{mn} \geq \left ( \frac{m^{n-1}}{ne} \right )^{m^n}
$$
by Stirling's formula.

For the upper bound, we note that one can think of ``building up'' a linear extension by choosing the next element at each step, and the set of possible ``next elements'' is an antichain.  It is well known that $[m]^n$ is Sperner, i.e., its largest antichain is its largest rank-level.  It is also rank-symmetric and unimodal; therefore, its largest rank-level is the ``middle'' one, i.e., rank $\floor{n(m-1)/2} + 1$.  This is the number of ways to write $\floor{n(m-1)/2} + 2$ as a sum of $n$ integers in $[m]$, a quantity bounded by
$$
\binom{\floor{n(m-1)/2} + 1}{n-1} \leq \frac{(nm/2)^{n-1}}{(n-1)!} \leq \left ( \frac{m e}{2} \right)^{n-1}.
$$
Therefore, an upper bound for the number of linear extensions is given by the $m^n$-th power of this quantity.
\end{proof}

In the case that the poset in question is not necessarily equilateral, it is not hard to extend the above argument in the case of an upper bound.

\begin{lemma} \label{notsquare} The number of linear extensions of $P = [a_1]\times \cdots \times [a_k]$ is at most
$$
\left ( \frac{|P|}{\max_j a_j} \right )^{|P|}.
$$
\end{lemma}
\begin{proof} We show, inductively, the claim that the largest antichain in $P$ has size at most $|P|/\max_j a_j$.  This clearly implies the lemma.  First, note that the $k=1$ case is trivial.  Now, suppose the claim is true for $k-1$.  Since $P$ is Sperner, we need only consider the size of its rank-levels.  The cardinality $|\cW_s|$ is given by the number of solutions to
$$
x_1 + \cdots + x_k = s + k - 1,
$$
where $1 \leq x_j \leq a_j$ is an integer.  Without loss of generality, we may assume that $a_k = \min_j a_j$.  Then we may write
$$
|\cW_s| = \sum_{x_k = 1}^{a_k}\left  |\left \{(x_1,\ldots,x_{k-1}) \in [a_1] \times \cdots \times [a_{k-1}] : \sum_{j=1}^{k-1} x_j = T - x_k \right \} \right|.
$$
Invoking the inductive hypothesis, we have
$$
|\cW_s| \leq a_k \cdot \frac{\prod_{j=1}^{k-1} a_j}{\max_{j=1}^{k-1} a_j} = \frac{|P|}{\max_{j=1}^k a_j}.
$$
\end{proof}
\bigskip

\section{Degree Sequence Extremes}

Think of a linear extension $L$ of a poset $P$ being ``built up'' an element at a time, so that the $j^\textrm{th}$ element of $L$ is ``added at time $j$''.  Then, given an integer $k$ with $1 \leq k \leq |P|$, define a {\it pit} to be any minimal element of
$$
P \setminus L^{[k]},
$$
i.e., a pit is an element of $P$ which could be ``next'' in any linear extension which agrees with $L$ up to time $k$.  Let $\pits(L,k)$ denote the number of pits at time $k$.

Note that the map $L^\prime_j : [m]^n \rightarrow [m^n]$ which agrees with $L$ except on the $j^\textrm{th}$ and $(j+1)^\textrm{st}$ elements, but has those two elements swapped, is a {\it bona fide} extension of $P = [m]^n$ if and only if $L^{j+1} \nsucc L^j$.  Therefore, $\deg(L)$ is the number of $j$ so that $L^{j+1} \nsucc L^j$, i.e., the number of ``good'' times $j$.  If the $(j+1)^\textrm{st}$ element of $L$ is chosen uniformly at random, then the probability that $j$ is bad is $d^+(L^j)/\pits(L,j)$, where $d^+(\bv)$ is the number of elements covering $\bv$ in the poset.  Note that good times correspond exactly to jumps.

\begin{theorem} $\max_{L \in \cT([m]^n)} \deg(L) = m^n - 3$.
\end{theorem}

\begin{proof} We begin by noting that there are $m^n-1$ possible times $j$.  At least two are always bad: the first and the last.  Indeed, the first is bad because $L^1 = (1,\ldots,1)$ and $L^2$ must cover this.  If $L^2 \nsucc L^1$, then there is a $z \in [m]^n$ with $L^1 < z < L^2$ so that $L(z) > 2$, contradicting that fact that $L$ is a linear extension of $[m]^n$.  Similarly, if $L^{m^n} \nsucc L^{m^n - 1}$, then there is a $z \in [m]^n$ with $L^{m^n - 1} < z < L^{m^n}$ so that $L(z) < m^n - 1$, again, a contradiction.  Therefore, $\deg(L) \leq m^n-3$ for all $L$.

On the other hand, this bound is achievable.  Indeed, let $L$ be the linear extension so that $r(x) < r(y)$ implies $L(x) < L(y)$, with the lexicographic ordering within each $\cW_s$.  Within each rank-level, there are no bad times.  The only possible bad times occur during the transition from $\cW_s$ to $\cW_{s+1}$.  However, note the first element of $\cW_s$ is the vector $(x_1,\ldots,x_n)$ constructed inductively by letting
$$
x_k = \max \{s + n - 1 - \sum_{j=1}^{k-1} x_j - m(n-k),1\}
$$
for $k = 1,\ldots,n$, since this is the smallest quantity $\geq 1$ so that writing $x_{k+1} + \cdots + x_n = s + n - 1 - \sum_{j=1}^{k} x_j$ does not require any integers greater than $m$ on the left-hand side.  The first coordinate will then be $s - (n-1)(m-1)$ if $s > (n-1)(m-1)$ and $1$ otherwise.  Similarly, the last element of $\cW_s$ is the vector $(x_1,\ldots,x_n)$ constructed inductively by letting
$$
x_k = \min \{s + k - 1 - \sum_{j=1}^{k-1} x_j,m\}
$$
for $k = 1,\ldots,n$, since this is the largest quantity $\leq m$ so that writing $x_{k+1} + \cdots + x_n = s + n - 1 - \sum_{j=1}^{k} x_j$ does not require any integers less than $1$ on the left-hand side.  The first coordinate will then be $s$ if $s \leq m$ and $m$ otherwise.

The upshot is that in level $s$, the last element has first coordinate $s$ if $s \leq m$ and $m$ otherwise, and then the next element in the ordering has first coordinate $s + 1 - (n-1)(m-1)$ if $s + 1> (n-1)(m-1)$ and $1$ otherwise.  But 
$$
\min \{s,m\} > \max \{s + 1 - (n-1)(m-1), 1\}
$$
whenever $s > 1$, $m > 1$, $(n-1)(m-1) > 1$, and $m > s + 1 - (n-1)(m-1)$.  The first two inequalities are always satisfied, except for the first time $j = 1$.  The third inequality only fails when $n = m = 2$, where it is easy to check that the claimed upper bound is correct.  The fourth inequality is satisfied whenever $s < n(m-1)$, i.e., $s$ is not the second-to-largest or largest rank.  Therefore, no element is succeeded by an element that covers it, except for the first and the second-to-last elements in the ordering.  This implies that $\deg(L) = m^n - 3$.
\end{proof}

For the minimum degree, we appeal to \cite{J95}.

\begin{theorem}[Jung '95] $s([m]^n) = \min_{L \in \cT([m]^n)} \deg(L) = m^{n-1} -1.$
\end{theorem}
\bigskip

\section{Average Jump Number}

Let $\BH[X]$ denote the entropy of the random variable $X$, $\BE[X]$ the expectation of $X$, and $\BP[A]$ the probability of the event $A$.  We will write $\lg$ for the logarithm base 2.

The following theorem shows that, in fact, $\cT([m]^n)$ is ``almost regular'' whenever $n = o(\log m)$.

\begin{theorem} The average degree in the graph $\cT([m]^n)$, i.e., $\bar{s}([m]^n)$, is at least
$$
m^n \left (1 - \sqrt{\frac{48 \lg n}{\lg m}} \right ),
$$
for sufficiently large $m$ and $n$.
\end{theorem}

\begin{remark} Here we can take ``sufficiently large'' to mean that $n \geq 2$ and $m \geq 41$.
\end{remark}

\begin{proof} If we let $\bL$ be a linear extension of $[m]^n$ chosen uniformly at random, then
\begin{equation} \label{eq2}
\BH[\bL] = \lg (A(m,n)) \geq (n-1) m^n \left (\lg (m) - \beta(n) \right).
\end{equation}
where $\beta(n) = (1+\lg n)/(n-1)$.  Suppose that, for all $k \in \cI \subset [m^n]$, we have
$$
\BH[\bL^{k+1} | \bL^{[k]}] < (n-1) \lg m - Q.
$$
Then
\begin{align*}
\sum_{k=1}^{m^n-1} \BH[\bL^{k+1} | \bL^{[k]}] & < |\cI| ((n-1) \lg m - Q) + (m^n - |\cI|) \lg T\\
& \leq S ((n-1) \lg m - Q)+ (m^n - |\cI|) (n-1) \lg m \\
& \leq (n-1) m^n \lg m  - Q |\cI|
\end{align*}
where $T = \left (\frac{me}{2}\right)^{n-1}$ is an upper bound for the size of the largest antichain in $[m]^n$.  Hence,
$$
(n-1) m^n ( \lg m - \beta(n) ) \! \leq \BH[\bL] = \!\!\sum_{k} \BH[\bL^{k+1} | \bL^{[k]}] \!< (n-1) m^n \lg m  - Q |\cI|,
$$
that is,
\begin{align*}
Q|\cI| & < (n-1) m^n \lg m - (n-1) m^n ( \lg m - \beta(n) ) \\
& = \beta(n) (n-1) m^n = m^n (1+ \lg n).
\end{align*}
Therefore, the number of times $k$ when $\BH[\bL^{k+1} | \bL^{[k]}] \geq (n-1) \lg m - Q$ is at least $m^n (1 - (1+ \lg n))/Q)$.

Now, suppose that, for some value $x$ and set $S \subset [m]^n$, the variable $\bL^{k+1} | (\bL^{[k]}= S)$ takes on the value $x$ with probability more than $\epsilon$.  Then
$$
\BH[\bL^{k+1} | \bL^{[k]} = S] \leq -\epsilon \lg \epsilon + (1 - \epsilon) (n-1) \lg m.
$$
If this event occurs (with respect to $\bL^{[k]}$) with probability at least $1- \epsilon$, where $\epsilon \leq 1/2$, then
\begin{align*}
\BH[\bL^{k+1} | \bL^{[k]}] & \leq (1-\epsilon) (-\epsilon \lg \epsilon + (1 - \epsilon) (n-1) \lg m) + \epsilon (n-1) \lg m \\
& = (n-1) \lg m  - \epsilon (1-\epsilon) (\lg \epsilon + (n-1) \lg m) \\
& \leq (n-1) \lg m  - \epsilon (n-1) \lg m / 2.
\end{align*}
Hence, if $\BH[\bL^{k+1} | \bL^{[k]}] \geq (n-1) \lg m - Q$, then $\epsilon \leq 2Q/((n-1) \lg m)$.  Then for at least $m^n (1 - (1+\lg n)/Q)$ times $k$, with probability at least $1 - 2Q/((n-1) \lg m)$ for $\bL^{[k]}$, no outcome of $\bL^{k+1} | \bL^{[k]}$ is likelier than probability $2Q/((n-1) \lg m)$.  (Note that, by the choice of $Q$ below, $\epsilon \leq 1/2$.)

Now, one can express the expected number of neighbors of the vertex $\bL$ in $\cT([m]^n)$ as
\begin{align*}
\BE[\deg(\bL)] &= \sum_{k=1}^{m^n - 1} \BP[\bL^{k+1} \nsucc \bL^k] \\
& \geq m^n \left (1 - \frac{1+\lg n}{Q} \right) \!\! \left (1 - \frac{2Q}{(n-1) \lg m} \right ) \!\! \left (1 - \frac{2 n Q}{(n-1) \lg m}\right ) \\
\end{align*}
where the third factor has an extra factor of $n$ because no element of $[m]^n$ is covered by more than $n$ elements. Continuing the computation, we find that
\begin{align*}
\BE[\deg(\bL)] & \geq m^n \left (1 - \frac{1+\lg n}{Q} - \frac{2(n+1) Q}{(n-1) \lg m} \right ) \\
& \geq m^n \left (1 - \frac{2\lg n}{Q} - \frac{6 Q}{\lg m} \right ).
\end{align*}
Taking $Q = \sqrt{\lg n \lg m/3}$ yields
$$
\BE[\deg(\bL)] \geq m^n \left (1 - \sqrt{\frac{48 \lg n}{\lg m}} \right ).
$$
\end{proof}

We can apply the fact that the maximum degree in $\cT([m]^n)$ is $m^n - 3$ to get the following more quantitative result.

\begin{cor} For all but at most a fraction $(48 \lg n/\lg m)^{1/4}$ vertices $L$ of $\cT([m]^n)$, we have
$$
\deg(L) \geq m^n \left (1 - \left ( \frac{48 \lg n}{\lg m} \right )^{1/4} \right).
$$
In particular, if $n = e^{o(\log m)}$, then the graph $\cT([m]^n)$ is almost regular.
\end{cor}

Of course, one can use this statement to conclude that the set of pits at most times in a random linear extension is, with high probability, a constant fraction of the size of a maximal antichain.  The conclusion is fairly weak, however, and only applies in the $\log n = o(\log m)$ regime.  On the other hand, a simpler proof yields a better conclusion, as the next result demonstrates.

\begin{prop} \label{ple} If a linear extension $L$ of $[m]^n$ is chosen uniformly at random, the expected fraction of times when $\pits(L,k) < 2^{-R} \left (\frac{me}{2}\right)^{n-1}$ is at most $(1+ \lg n)/R$.
\end{prop}
\begin{proof}
Note that
\begin{align*}
\BH[\bL] & = \sum_{k=1}^{m^n-1} \BH[\bL^{[k+1]} | \bL^{[k]}] \\
& = \sum_{k=1}^{m^n-1} \BH[\bL^{k+1} | \bL^{[k]}] \\
& \leq \sum_{k=1}^{m^n-1} \BE [\lg(\pits(\bL,k))] \\
& = \sum_{k=1}^{m^n-1} \sum_{r=1}^{T} \lg{r} \cdot \BP[\pits(\bL,k)=r].
\end{align*}
Hence,
\begin{align*}
\BH[\bL] & \leq \sum_{r=1}^{T} \lg{r} \sum_{k=1}^{m^n-1} \BP[\pits(\bL,k)=r] \\
& = \sum_{r=1}^{T} \lg{r} \cdot \BE[\# k \textrm{ such that } \pits(\bL,k)=r]
\end{align*}
Combining this with (\ref{eq2}) and writing $\eta(r) =\BE[\# k \backepsilon \pits(\bL,k)=r]$, we find
\begin{equation} \label{eq3}
\sum_{r=1}^{T} \eta(r) \lg{r} \geq (n-1) m^n ( \lg m - \beta(n) ).
\end{equation}
Suppose that $\sum_{r=\epsilon T}^T \eta(r) < m^n (1 + (1+\lg n)/\lg \epsilon)$ for some $\epsilon \in (0,1)$.  Then
\begin{align*}
\sum_{r=1}^{T} \eta(r) \lg{r} &= \sum_{r=\epsilon T}^{T} \eta(r) \lg{r} + \sum_{r=1}^{\epsilon T} \eta(r) \lg{r} \\
& \leq \lg{T} \sum_{r=\epsilon T}^{T} \eta(r) + (m^n - \sum_{r=\epsilon T}^{T} \eta(r)) \lg{\epsilon T} \\
& < - \lg{\epsilon} \left ( m^n + \frac{(1+ \lg n) m^n}{\lg \epsilon} \right ) + m^n \lg{\epsilon T} \\
& = m^n \lg T - (1+\lg n) m^n\\
& = (n-1) m^n (\lg m - \beta(n)),
\end{align*}
contradicting (\ref{eq3}).  Therefore, the expected fraction of times $k$ when $\pits(\bL,k) < \epsilon T$ is at most $- (1+ \lg n)/\lg \epsilon$.
\end{proof}

By applying Markov's inequality we immediately conclude the following.

\begin{cor} For $\Delta \geq 1$, the probability that $\bL$ has more than 
$$
\frac{\Delta (1+ \lg n)}{R}\cdot |P|
$$
times $k$ when $\pits(\bL,k) < 2^{-R} \left (\frac{me}{2}\right)^{n-1}$ is at most $\Delta^{-1}$.
\end{cor}

Note that this statement has force even when $n = e^{\Omega(\log m)}$.  For example, if $R = n = \sqrt{m}$ and $\Delta = m^{1/4}/(1+ \lg n)$, we conclude that the probability that $\bL$ has more than $m^{-1/4} |P|$ times $k$ when $\pits(\bL,k) < .7^{\sqrt{m}} |P|$ is $\tilde{O}(m^{-1/4})$.
\bigskip

\section{Notes}

We conjecture that, even when the dimension is large, the average jump number is close to the maximum jump number.

\begin{conj} $\bar{s}([m]^n) = m^n (1-o_{m,n}(1))$ for all $m$, $n$ with $m^n \rightarrow \infty$.
\end{conj}

It would be interesting to describe more precisely the integer sequence $\pits(L,k)$ for random choices of $L$.  Proposition \ref{ple} gets at this issue, but there is obviously much left to ask.  Is the sequence $\pits(L,k)$ unimodal with high probability after convolution with a short interval in $k$?  What are the higher moments of random jump numbers?  To what extent is the statistic $\bar{s}(\cdot)$ respected by Cartesion product?  Which posets have their average jump number substantially different from the maximum jump number?  What does it say about the relationship between the height and width of a poset if one knows that $\bar{s}(P)$ is close to the maximum?
\bigskip

\section{Acknowledgments}

Thank you to Angelika Steger and ETH-Z\"{u}rich for their hospitality and support for the months during which this work was done.  Thank you also to Jochem Giesen and Eva Schuberth for providing valuable insights.  Finally, gratitude to Graham Brightwell for allowing me to pick his encyclopedic knowledge of poset theory.

\end{document}